\documentclass[11pt]{amsart}

\newcounter{mtheorem}
\newtheorem{mtheorem}[mtheorem]{Theorem}

\newtheorem{theorem}{Theorem}[section]
\newtheorem{lemma}[theorem]{Lemma}
\newtheorem{prop}[theorem]{Proposition}
\newtheorem{corollary}[theorem]{Corollary}

\theoremstyle{definition}
\newtheorem{definition}[theorem]{Definition}
\newtheorem{example}[theorem]{Example}
\newtheorem{xca}[theorem]{Exercise}

\theoremstyle{remark}
\newtheorem{remark}[theorem]{Remark}

\numberwithin{equation}{section}

\newcommand{\abs}[1]{\lvert#1\rvert}

\newcommand{\blankbox}[2]{%
  \parbox{\columnwidth}{\centering
    \setlength{\fboxsep}{0pt}%
    \fbox{\raisebox{0pt}[#2]{\hspace{#1}}}%
  }%
}
\newcommand{\beq}{\begin{equation}}
\newcommand{\eeq}{\end{equation}}
\newcommand{\bea}{\begin{eqnarray}}
\newcommand{\eea}{\end{eqnarray}}
\newcommand{\cy}{{C}alabi-{Y}au\ }
\newcommand{\ka}{{K}\"ahler\ }

\newcommand{\slg}{special {L}agrangian\ }

\newcommand{\sleg}{special {L}egendrian\ }
\newcommand{\tslg}{$\theta$-special Lagrangian\ }
\newcommand{\tsleg}{$\theta$-special {L}egendrian\ }
\newcommand{\tslegg}{$\theta$-special {L}egendrian}
\newcommand{\C}{\mathbb{C}}
\newcommand{\R}{\mathbb{R}}
\newcommand{\Q}{\mathbb{Q}}
\newcommand{\Z}{\mathbb{Z}}
\newcommand{\N}{\mathbb{N}}
\newcommand{\ra}{\rightarrow}

\newcommand{\lra}{\longrightarrow}
\newcommand{\la}{\lambda}
\newcommand{\no}{\noindent}
\newcommand{\sn}{\operatorname{sn}}
\newcommand{\Ke}{\operatorname{Ke}}
\newcommand{\cnd}{\operatorname{cn}}
\newcommand{\dn}{\operatorname{dn}}
\newcommand{\diag}{\operatorname{diag}}
\newcommand{\sech}{\operatorname{sech}}
\newcommand{\hess}{\operatorname{Hess}}

\begin{document}

\title{Special Lagrangian cones}

\author{Mark Haskins}
\address{Department of Mathematics, University of Texas, Austin TX 78712}
\email{mhaskins@math.utexas.edu}
\thanks{The author was supported in part by the States of Guernsey Education Council.
He would like to thank the Institute for Advanced Study and Pomona College for their
hospitality.}


\date{May 16, 2000.}



\begin{abstract}
We study homogeneous special Lagrangian cones in $\C^n$ with isolated singularities.
Our main result constructs an infinite family of special Lagrangian cones in $\C^3$ each
of which has a toroidal link. We obtain a detailed geometric description of these tori. 
We prove a regularity result for special Lagrangian cones in $\C^3$ with a spherical link --
any such cone must be a plane. We also construct a one-parameter family of asymptotically
conical special Lagrangian submanifolds from any special Lagrangian cone.
\end{abstract}

\maketitle

\section{Introduction}

Let $Y$ be \cy manifold of complex dimension $n$ with \ka form $\omega$ and non-zero
parallel holomorphic $n$-form $\Omega$ satisfying the normalization condition
$\omega^n/ n! = (-1)^{n(n-1)/2}(i/2)\Omega \wedge \bar{\Omega}$. 
Then $\textrm{Re\ }(\Omega)$ is a calibrated form, 
whose calibrated submanifolds are called special Lagrangian submanifolds \cite{HL:1982}.

Moduli spaces of \slg submanifolds (and of other calibrated submanifolds) have appeared
recently in string theory \cite{becker}, \cite{SYZ:1996}. 
On physical grounds, Strominger, Yau and Zaslow argued 
\cite{SYZ:1996} that a \cy manifold $Y$ with a mirror partner $\hat{Y}$ admits a (singular) 
fibration by \slg tori, and that $\hat{Y}$ should be obtained by compactifying the dual fibration.
To make this idea rigorous one needs to have control over the singularities and compactness
properties of families of \slg submanifolds. In dimensions three and higher these properties are
not well understood. 

Motivated by these problems we study the simplest isolated singularities of \slg
varieties -- homogeneous cones in $\C^n$ with an isolated singularity. These
are also local models for more general singularities in that they are possible
tangent cones to \slg currents at singular points. 

We introduce the notion of a \tsleg submanifold -- a special class of minimal $(n-1)$-dimensional submanifolds -- of $S^{2n-1}(1)$, and characterize \tslg cones in $\C^n$ as those cones $C$ whose links
$L = C \cap S^{2n-1}$ are \tsleg submanifolds of $S^{2n-1}$ (Proposition \ref{tsleg}). 
From any \sleg link, in addition to a \slg cone, we obtain a one-parameter family of asymptotically conical \slg (possibly immersed) submanifolds.

\begin{mtheorem}
\label{ac_slg}
Let $\Sigma^{n-1}$ be a \tsleg submanifold of $S^{2n-1}(1) \subset \C^n$. 
Let $\Sigma_{d}$ ($d\in \R $) denote the set 
$$ \left \{ (zp \in \C^n): p \in \Sigma, z\in \C, \textrm{where}\ {\rm Im}(z^n)=d, \ \arg{z}\in [0,\frac{\pi}{n}] \right \}.$$ 
Then 

\no
(i) $\Sigma_d$ is a \tslg variety.

\no
(ii) $\Sigma_0 = C(\Sigma) \cup C(e^{i\pi/n}\Sigma)$ where $C(\Sigma)$ denotes the cone on $\Sigma$. 

\no
(iii) $\Sigma_d$ is asymptotically conical, with two ends $\Sigma$ and $e^{i\pi/n} \Sigma$.
\end{mtheorem}

In the case of \slg cones in $\C^3$, results  
of Yau \cite{yau} and others put restrictions on three-dimensional \slg cones.
For example, we obtain:

\begin{mtheorem}
\label{link=sphere}
Let $C$ be a homogeneous \slg cone in $\C^3$, with $L=C \cap S^5(1)$ a (possibly immersed) sphere.
Then $C$ must be a \slg plane.
\end{mtheorem}

A simple corollary of this theorem is a regularity result for homogeneous
solutions of the \slg graph equation in dimension three. The theorem is
also sharp in the following two senses. The analogous statement in $\C^4$ is
false, as recent examples of Chen \textit{et al.} \cite{chen} demonstrate. Moreover, if the link type is
a torus not a sphere then even in dimension three there are nontrivial \slg cones.
Our main result gives an abundance of such cones.

\begin{mtheorem}
\label{main_cone_thm}
There exists a countably infinite family of non-isometric \slg cones in $\C^3$. 
Each cone has link an embedded torus which is invariant under some $S^1 \subset \textrm{SU(3)}$.
\end{mtheorem}

Each special Lagrangian cone in $\C^3$ also gives rise to
other calibrated cones (Lemma \ref{cal_cones}). For example, from each special Lagrangian 
cone in $\C^3$ with an isolated singularity we associate: an associative cone in 
$\R^7$ with an isolated singularity, a coassociative cone in $\R^7$ and a Cayley
cone in $\C^4$ with a line of singularities and special Lagrangian cones in $\C^{n+3}$
with singularities along a real $n$-plane.

The strategy for the construction of the special Lagrangian cones in $\C^3$ is as follows.
By exploiting the connection between harmonic maps and minimal surfaces we
construct a two-parameter family $u_{\alpha,J}$ of special Legendrian immersions
$\R^2 \ra S^5(1)$.
Harmonic maps from two dimensional domains to Lie groups and 
symmetric spaces have a rich structure, with relations to infinite dimensional
completely integrable systems and loop groups \cite{Gu:1997}, \cite{uhlen:chiral}.
In the case of $S^1$-equivariant
harmonic maps to spheres, a finite dimensional completely integrable system -- the
C. Neumann system -- appears. Several geometric features of the harmonic map have
nice interpretations in terms of conserved quantities of this system. The mechanical interpretation
of the Legendrian condition is not so clear, but nonetheless we are able to
obtain minimal Legendrian immersions from certain solutions of the Neumann system.

\begin{mtheorem}
\label{min_lag_family}
For each $\theta\in [0,2\pi)$ there exists a $2$-parameter family $u_{\alpha,J}$, 
 $(\alpha,J) \in [0,1] \times [0,1/3\sqrt{3}]$, of \tsleg immersions $\R^2 \ra S^5(1)$
with the following properties:\\
(i) The immersion $u_{\alpha,J}$ is invariant under the $1$-parameter subgroup of $\textrm{SU(3)}$ generated by $A=\diag{i(1,\alpha, -1-\alpha)}\in su(3)$. \\
(ii) For $\alpha = 1$ and $J = 1/3\sqrt{3}$ these immersions
all describe Clifford tori, but otherwise all the immersions are geometrically distinct.\\
(iii) The family $u_{\alpha,J}$ contains all \tsleg Legendrian immersions of the form 
\eqref{equi} which cover a torus.
\end{mtheorem}

To obtain tori from these immersions, 
in Proposition \ref{period_prop} we examine the conditions
under which $u_{\alpha,J}$ is doubly periodic with respect to some lattice.
By examining the special cases $u_{J,0}$ and $u_{0,\alpha}$ we deduce
\begin{mtheorem}
\label{tori}
(i) For $\alpha \in \Q \cap (0,1]$, the immersion $u_{0,\alpha}$ is doubly periodic and
hence gives rise to a minimal Legendrian torus.\\
(ii) For a dense set of $J \in (0,1/3\sqrt{3})$ the immersion $u_{J,0}$ is doubly periodic and
hence gives rise to a minimal Legendrian torus. 
\end{mtheorem}

For the family $u_{0,\alpha}$, referred to in part (i) of the previous theorem, we give detailed information 
about the geometry (e.g. conformal structure, maximum and minimum values of the Gauss curvature and embeddedness) of the corresponding surfaces. As a corollary we find (Theorem \ref{almost_flat}) 
that there are embedded \lq almost flat\rq\  minimal Legendrian tori. These tori demonstrate sharpness
of two pinching results of Yau on minimal Lagrangian (Legendrian) immersions into $\C P^2$ ($S^5$).

The paper is organized as follows. In Section 2 we recall basic facts about \slg geometry 
in $\C^n$, introduce the notion of \sleg in $S^{2n-1}$ and characterize \slg cones in
terms of \sleg links. In Section 3 we recall basic facts from harmonic map theory: principally
the relation with minimal surfaces and the appearance of the C. Neumann system in
$S^1$-equivariant harmonic maps into spheres. In Section 4 we study which solutions of the Neumann system give rise to special Legendrian immersions, give explicit parametrisations of these solutions and study the geometry of these immersions. In Section 5 we study the periodicity conditions for these immersions and  hence are able to deduce our main results.

\section{\slg cones in $\C^n$}
\subsection{Special Lagrangian geometry in $\C^n$}
Special Lagrangian geometry is an example of a \textit{calibrated geometry} \cite{HL:1982}.
We review some elementary facts about calibrations and special Lagrangian geometries in $\C^n$
in particular (see \cite{HL:1982} for further details).

Each calibrated geometry is a distinguished class of minimal submanifolds of a Riemannian manifold $(M,g)$ associated with a closed differential $p$-form $\phi$ of comass one.

For each $m\in M$, the \textit{comass} of $\phi$ is defined to be
\begin{displaymath}
\|\phi\|_m^* = \sup\{ <\phi_m,\xi_m> : \xi_m \ \textrm{is a unit simple } p\textrm{-vector at m}\}.
\end{displaymath}
\no
In other words, $\|\phi\|_m^{*}$ is the supremum of $\phi$ restricted to the Grassman of
oriented $p$-dimensional planes $G(p,T_m M)$, regarded as a subset of $\Lambda^{p}T_mM$.

To any form of comass one there is a natural subset of $G(p,TM)$
\begin{displaymath}
G_m(\phi)  = \{\xi_m \in G(p,T_mM) : <\phi_m,\xi_m> = 1 \} ,
\end{displaymath}
that is, the collection of oriented $p$-planes on which $\phi$ assumes its maximum.
These planes are the planes \textit{calibrated by $\phi$}.
An oriented $p$-dimensional submanifold of $(M,g)$ is \textit{calibrated by $\phi$}
if its tangent plane at each point is calibrated.  

The key property of calibrated submanifolds is that they are \textit{homologically
volume minimizing}

\begin{lemma}[Harvey and Lawson \cite{HL:1982}]
Let $(M,g,\phi)$ be a calibrated geometry, and suppose $S$ is a calibrated submanifold
(possibly with boundary). Then for any oriented $p$-dimensional submanifold $\hat{S}$ homologous to $S$
$$ \textrm{vol}(S) \le \textrm{vol}(\hat{S})$$
with equality if and only if $\hat{S}$ is also calibrated (by $\phi$).
\end{lemma}
\no
Let $z_1, \ldots, z_n$ denote standard complex coordinates on $\C^n$. 
For any $\theta \in [0,2\pi)$ the real $n$-form 
$$\alpha_{\theta} = \textrm{Re}(e^{i\theta}dz^1 \wedge \ldots \wedge dz^n)$$ 
is a calibrated form, called the \textit{\tslg calibration} on $\C^n$.

For the proof that $\alpha_{\theta}$ has comass one see \cite{HL:1982}.
A \textit{\tslg plane} (we will sometimes abbreviate this as $\theta$-SLG) is an oriented $n$-plane calibrated by the form $\alpha_{\theta}$. 
A useful characterization of the \tslg planes is
\begin{lemma}
An oriented $n$-plane $\xi$ in $\C^n$ is \tslg (for the correct choice 
of orientation) if and only if
\begin{enumerate}
\item
$\xi$ is Lagrangian with respect to the standard symplectic form 
$\omega = \sum{dx^i \wedge dy^i}$, (i.e. $\omega$ restricts to zero on $\xi$) and
\item 
$\beta_{\theta} := \textrm{Im}(e^{i\theta}dz^1 \wedge \ldots \wedge dz^n)$ restricts to zero on $\xi$.
\end{enumerate}
\end{lemma}
One reason for considering the whole $S^1$-family of \slg calibrations is the following 
result of Harvey and Lawson (Proposition 2.17 of \cite{HL:1982}):

\begin{prop}
\label{mlag}
A connected oriented Lagrangian submanifold $S \subset \C^n$ is minimal (i.e. it is
a critical point of volume, or its mean curvature $H$ vanishes) if and only if
$S$ is \tslg for some $\theta$.
\end{prop}

\subsection{Regular cones and \sleg links}
For any compact connected oriented embedded submanifold
$\Sigma \subset S^{n-1}(1)\subset \R^n$ define the \textit{cone on $\Sigma$}, 
$$ C(\Sigma) = \{ tx: t\in \R^{\ge 0}, x \in \Sigma \}.$$
A cone $C$ in $\R^n$ is $\textit{regular}$ if there exists $\Sigma$ as above so
that $C=C(\Sigma)$, in which case we call $\Sigma$ the \textit{link} of the cone $C$.
$C(\Sigma) - 0$ is an embedded smooth submanifold, but $C(\Sigma)$ has an isolated
singularity at $0$ unless $\Sigma$ is a totally geodesic sphere.

To characterize the links of regular \slg cones we need to introduce some
geometric structures on the unit sphere $S^{2n-1}$ in $\C^n$. As a convex hypersurface
in a K\"ahler  manifold \cite{mcd}, $S^{2n-1}(1)$ inherits a 
\textit{contact form}, that is, a $1$-form $\gamma$ so that
\begin{equation}
\label{contact}
\gamma \wedge d\gamma^{n-1} \neq 0.
\end{equation}
Let $X$ denote the Euler vector field $x\cdot\partial /\partial x$ on $\C^n$ and $\omega$ denote the 
standard symplectic form on $\C^n$. Then the contact form on $S^{2n-1}(1)$ is
$$ \gamma = \iota_X\omega | _{S^{2n-1}}.$$
Associated with $\gamma$ is the \textit{contact distribution}, the hyperplane
field $\ker{\gamma} \subset TS^{2n-1}$. The condition (\ref{contact}) on $\gamma$
ensures that the distribution $\ker{\gamma}$ is not integrable. The maximal dimensional
integral submanifolds (i.e. submanifolds on which $\gamma$ restricts to zero) of the distribution
are $(n-1)$-dimensional and are called \textit{Legendrian submanifolds}.

The relevance of Legendrian submanifolds of the sphere can be see from the next result
whose proof is standard.

\begin{lemma}
Let $\Sigma$ be an $(n-1)$-dimensional submanifold of $S^{2n-1}(1)$. Then 
$C(\Sigma)$ is Lagrangian if and only if $\Sigma$ is Legendrian.
\end{lemma}

For any $p$-form $\phi$ on $\R^n$ define the \textit{normal part} of $\phi$ by
$$ \phi_N = \iota_X \phi,$$ 
where $X$ again denotes the Euler vector field on $\C^n$.
In particular, $\alpha_{\theta,N}$ denotes the normal
part of the \tslg calibration $\alpha_{\theta}$.

An oriented $(n-1)$-dimensional submanifold $\Sigma$ of $S^{2n-1}(1)$ is a \textit{\tsleg submanifold}
if at each point of $\Sigma$, $\alpha_{\theta,N}$ restricts to the volume form on $\Sigma$.

\begin{prop}
\label{tsleg}
A regular cone $C=C(\Sigma)$ in $\C^n$ is \tslg if and only if $\Sigma$ is \tslegg.
\end{prop}
\begin{proof}
This is essentially a special case of Theorem 5.6 of \cite{HL:1982}. We sketch the proof.
For any constant $p$-form $\phi$ on $\R^n$ define the \textit{tangential part} of $\phi$ to be
$$ \phi_T = \iota_X\left(\frac{x}{|x|}\cdot dx \wedge \phi\right).$$ Then $\phi$ decomposes as
\begin{equation} 
\phi = \phi_T + \frac{x}{|x|}\cdot dx \wedge \phi_N
\end{equation}
where $\phi_N$ is the normal part of $\phi$ defined previously. When $d\phi=0$, 
restricting to the unit sphere it follows that
$$ d\phi_T = 0 \quad \textrm{and} \quad d\phi_N = p \phi_T.$$
Since $\phi_N(\xi) = \phi(x\wedge \xi)$ and $||\xi||= \|x\wedge \xi\|$
for any simple $(p-1)$-vector in $\Lambda^{p-1}x^{\perp}$, $\phi_N$
still has comass one (but since $\phi_N$ is not closed it is not a calibration itself).
Moreover,  submanifolds $\Sigma$ of $S^{n-1}(1)$ on which $\phi_N$ restricts to the volume form
are exactly those for which $C(\Sigma)$ is calibrated by $\phi$. Hence the result follows
by taking $\phi$ to be any of the \tslg calibrations $\alpha_{\theta}$.
\end{proof}

We also have the following Legendrian analogue of Proposition \ref{mlag}
\begin{prop}
\label{minlag_tsleg}
A connected oriented Legendrian submanifold of $S^{2n-1}(1)$ is minimal if and only if
it is \tsleg for some $\theta$.
\end{prop}

\begin{proof}
Let $\Sigma$ be a minimal Legendrian submanifold of $S^{2n-1}(1)$. It is a standard fact \cite{S:1968} that $C(\Sigma)$ is minimal if and only if $\Sigma$
is minimal in the unit sphere. Thus $C(\Sigma)$ is a minimal Lagrangian cone which from 
Proposition \ref{mlag} must be \tslg for some $\theta$. By the previous proposition this
implies $\Sigma$ is \tslegg. The converse is similar.
\end{proof}

We finish the section by proving Theorem \ref{ac_slg}, which gives a one-parameter family of asymptotically conical \slg varieties modeled on any \slg cone.
This result generalizes Theorem 3.5 of \cite{HL:1982} which is our
result in the special case  
$\Sigma = \left\{(x_1, \ldots ,  x_n) \in \C^n: x_i \in \R \ \textrm{with}\  \sum{x_i^2} =1\right\}$.

\begin{proof}[Proof of Theorem \ref{ac_slg}]
(i) By Proposition \ref{tsleg}, $C(\Sigma)$ is a $\theta$-SLG cone. 
By rotating $\Sigma$ by $A = \textrm{diag}(e^{-i\theta /n}, \ldots , e^{-i\theta /n})$ we can assume
$C(\Sigma)$ is $0$-SLG. Thus $\beta |_{C(\Sigma)} = 0$. Let $\phi: \Sigma \ra S^{2n-1}$ denote the inclusion
of $\Sigma$ in the sphere, and let $x_1, \ldots, x_{n-1}$ be local coordinates on $\Sigma$. Then
$\beta|_{C(\Sigma)} = 0$ is equivalent to
\begin{equation}
  \label{im_det}
  \textrm{Im}\left(\det{_{\C}(\phi, \frac{\partial \phi}{\partial x_1}, \ldots, \frac{\partial \phi}{\partial x_{n-1}})}\right)=0.
\end{equation}
Let $\Phi:\R \times \Sigma \ra \C^n$ be given by $\Phi(t,x)=f(t)\phi(x)$ where $f:\R \ra \C$ is some nonconstant
smooth complex valued function. It is straightforward to check that any such $\Phi$ gives rise to a
Lagrangian immersion to $\C^n$. 

Let $x_0 = t$ and for $i=0, \ldots ,n-1$ denote $\partial \Phi /\partial x_i$ 
by $\Phi_i$. Then for $j=1, \ldots , n-1$
$$ \omega(\Phi_0,\Phi_j) = \omega(\dot{f}\phi,f\phi_j) = {\rm Re}(\bar{f}\dot{f})\ \omega(\phi,\phi_i) + {\rm Im}(\bar{f}\dot{f})<\phi,\phi_i>=0$$
where the first term vanishes because $\phi$ is Legendrian and the second because $|\phi|^2=1$. For
$j,k = 1, \ldots ,n-1$ we have
$$ \omega(\Phi_j,\Phi_k)=\omega(f\phi_j,f\phi_k)=|f|^2\omega(\phi_j,\phi_k)=0$$
and so $\Phi$ is Lagrangian as claimed.

Now we claim that $\Phi$ is $0$-SLG if and only if $f$ satisfies 
\begin{equation}
  \label{zn}
 {\rm Im}\ (f^n)=d 
\end{equation}
for some real constant $d$. 
To prove this it is enough to show that $\beta |_{\Phi} = 0$ holds if and only if (\ref{zn}) is satisfied.
But $\beta |_{\Phi}=0$ is equivalent to 
\begin{equation}
\label{beta=0}
  \textrm{Im}\left( \det{_{\C}(\Phi_0, \ldots, \Phi_{n-1})}\right)=0.
\end{equation}
Since $C(\Sigma)$ is $0$-SLG, we have 
$   \textrm{Im}\left( \det{_{\C}(\phi, \phi_1, \ldots, \phi_{n-1})}\right)=0$
and hence
\begin{eqnarray*}
   \textrm{Im}\left( \det{_{\C}(\Phi_0, \ldots, \Phi_{n-1})}\right) & = &\textrm{Im}\left( \det{_{\C}(\dot{f}\phi, f\phi_1, \ldots, f\phi_{n-1})}\right)\\
& = & \textrm{Im}\left( \dot{f}f^{n-1}\det{_{\C}(\phi, \phi_1, \ldots, \phi_{n-1})}\right)\\
& = & \textrm{Re}\left( \det{_{\C}(\phi, \phi_1, \ldots, \phi_{n-1})}\right) \times \textrm{Im}(\dot{f}f^{n-1}).
\end{eqnarray*}
Thus (\ref{beta=0}) holds if and only if
$$ \textrm{Im} (\dot{f}f^{n-1}) = \frac{1}{n}\left( \textrm{Im}\frac{d}{dt}f^n \right) = \frac{1}{n}\frac{d}{dt} \textrm{Im}(f^n)=0.$$
Hence $\Phi$ is $0$-SLG if and only if $ {\rm Im}\ (f^n)=d $ as claimed.

Parts (ii) and (iii) are straightforward to verify.
\end{proof}

\subsection{Minimal Legendrian surfaces}
In dimension two, any \slg cone must be a union of \slg planes (since its link
must be a union of Legendrian geodesics in $S^3$). In the first interesting
case, namely \slg cones in $\C^3$, restrictions on the geometry and topology of the
allowable links follow from the next result, essentially due to Yau.

\begin{theorem}
\label{yau}
\cite{yam,yau}
Let $\Sigma$ be a minimal Legendrian surface of  $S^5(1)$. Then:

\no(i) If $\Sigma$ has genus zero, $\Sigma$ is totally geodesic. 

\no(ii) If $\Sigma$ is a complete nonnegatively curved surface, $\Sigma$ is a totally geodesic sphere
or a flat torus.

\no(iii) If $\Sigma$ is complete nonpositively curved surface then $\Sigma$ is a flat torus.
\end{theorem}

This theorem is the Legendrian analogue of a result of Yau on minimal Lagrangian immersions
into K\"ahler surfaces of constant holomorphic sectional curvature (e.g. $\C P^2$ with the
Fubini-Study metric). In fact, using Reckziegel's observation \cite{reck} about the local correspondence between minimal Legendrian immersions into $S^5$ and minimal Lagrangian immersions into $\C P^2$, one can deduce the Legendrian result from the Lagrangian one.

As a corollary of part (i) of Theorem \ref{yau} we deduce Theorem \ref{link=sphere}.
Applying this theorem to the special case of \slg graphs we deduce
\begin{corollary}
\label{slg-graph}
Any homogeneous degree 1 solution $u$ of the 3-dimensional \slg graph equation
$$ \Delta u = \det{\hess{(u)}}$$ is a quadratic function.
\end{corollary}

Theorem \ref{link=sphere} is sharp in the following two senses. Firstly, in $\C^4$ the analogous result
is false as recent examples of Chen \textit{et al.} show \cite{chen}. 
Secondly, in $\C^3$ there
are nontrivial \slg cones with link type a torus, the simplest example of which is the cone on a generalized Clifford torus. 

Let $T$ be the Lagrangian product $3$-torus contained in $S^5(1)$
$$ T = \{ z\in \C^3: |z_i|^2 = 1/3, \quad i=1,2,3 \},$$
and $T_{\theta}$ be the $2$-torus
$$ T_{\theta} = \{z\in T: \sum{\arg{z_i}} = \theta\}.$$
The $T_{\theta}$, the \textit{generalized Clifford tori}, are all flat minimal Legendrian tori and
$T_{\theta /3}$, $T_{\pi + \theta /3}$ are \tslegg. Harvey and Lawson discovered the cones on these tori
in a family of \slg level sets invariant under the maximal torus $T^2 \subset \textrm{SU(3)}$. 
This high degree of symmetry allowed them to explicitly write down solutions. In the next three sections
we shall find a family of nonisometric minimal Legendrian tori in $S^5(1)$, which
are invariant under an $S^1 \subset \textrm{SU(3)}$. This symmetry is still enough to allow us to give
quite explicit descriptions of these tori.

Special Lagrangian cones in dimension three with isolated singularities,
also give rise naturally to several 
related singular calibrated varieties. For example, remarks of Harvey-Lawson 
\cite{HL:1982} (IV.2.C. Remark 2.12) and Donaldson-Thomas \cite{donaldson} show
the following:
\begin{lemma}
\label{cal_cones}
If $X^3$ is a $0$-\slg variety, then:

\no
(i) $X \times \{pt\} \subset \C^3 \times \R $ is an associative variety \\
\no
(ii) $X \times \R \subset \C^3 \times \R$ is a coassociative variety \\
\no
(iii) $X \times \R \subset \C^3 \times \C$ is a Cayley variety.
\end{lemma}
\no
In cases (ii) and (iii) starting with a \slg cone with an isolated singularity we obtain \textit{cylindrical} cones, which have
a whole line of singularities. One also gets \slg cylindrical cones in $\C^{n+3}$
by taking the Cartesian product of a $3$-dimensional cone in $\C^3$ with a real $n$-plane
in $\C^n$.

\section{Harmonic Maps, Minimal Surfaces and the Neumann System}

We shall construct $S^1$-invariant minimal Legendrian tori in $S^5(1)$ by exploiting 
two relationships. The first is the connection between harmonic maps and minimal surfaces. 
The second is the link between $S^1$-equivariant harmonic maps into spheres and 
the C. Neumann system describing motion on the sphere under a quadratic potential.

\subsection{Harmonic Maps}
We recall some definitions and basic facts from harmonic map theory.
Suppose $M$ and $N$ are Riemannian manifolds. For any $C^1$ map $u: M \lra N$ define a smooth 
function $e(u)$, the \textit{energy density} of $u$, by $e(u)(x) = Tr(du_x^2)$. Define a functional
on $C^1(M,N)$, the \textit{total energy}, by $E(u)=\int_{M}e(u)\mu_M$, where $\mu_M$ is the 
Riemannian volume element of $M$. Critical points of $E$ are \textit{harmonic maps} from $M$ into $N$.



If $N$ is isometrically embedded in $\R^K$ then we can view a function $u:M \ra N$ as a function 
$u = (u^1,\ldots,u^K)$ into $\R^K$ with the constraint that $u(x)\in N$ for all $x\in M$. Then
\begin{equation}
  E(u) = \sum_{i=1}^{K} \int_M |\nabla u^i|^2 \mu_M.
\end{equation}
Extremals of $E$ subject to the constraint that $u(M)\subset N$ give us the harmonic maps to $N$. 
From this we see that the harmonic map equations can be written simply as
\begin{equation*}
\Delta u(x) \perp T_{u(x)}N, \quad u(x)\in N, \quad \forall x\in M.   
\end{equation*}
In the case that $N=S^n(1) \subseteq \R^{n+1}$ (with the metric induced by this inclusion)
this implies $\Delta u= \la u$ for some function
$\la$ on $M$. Taking the inner product of both sides with $u$ and using the constraint equation
$|u|^2 = 1$ we determine that $\la = (u,\Delta u)=-|du|^2$. Summarizing we have

\begin{lemma}
A smooth map $u:M\ra S^n(1)\subset \R^{n+1}$ is harmonic iff and only if
$u$ satisfies the equation
\begin{equation}
  \label{hsphere}
  \Delta u = - |du|^2 u.
\end{equation}
\end{lemma}
Finally we recall what happens to the harmonic map equations when we make a conformal change
of metric on the domain $M$. If $\widetilde{g}=\la^2 g$ then 
$\widetilde{g}^{-1}=\la^{-2}g^{-1}$, and $\widetilde{\mu}_M = \la^{m} \mu_M$. 
Hence $E_{\tilde{g}}(u)= \la^{m-2}E_g(u)$ and we see that $E$ is conformally invariant 
if and only if dim $M=2$. Therefore in dimension 2 harmonicity depends only on the structure of $M$ as a 
Riemann surface. In particular there is a natural quadratic differential $\Phi$, the
\textit{Hopf differential}. If $z$ is a local complex coordinate on $M$ then $\Phi = \phi(z)dz^2$
where
\begin{equation}
  \phi(z) = (u_z, u_z) = \frac{1}{4} \left( |u_x|^2 - |u_y|^2 - 2i (u_x,u_y) \right).
\end{equation}
Harmonicity of $u$ implies that $\Phi$ is holomorphic. If the Hopf differential vanishes
the map $u$ is \textit{conformal}. Moreover, we have the following connection with minimal surfaces:
\begin{prop}(\cite{ratto})
  $u$ is harmonic and conformal if and only if $u$ is a (branched) minimal immersion.
\end{prop}

\bigskip

\subsection{Equivariant Harmonic Maps and the Neumann System}

For harmonic maps from $\R^2$ to $S^5(1)\subseteq \C^3$ 
(where both $\R^2$ and $S^5(1)$ are given their standard metrics) of the special form
\begin{equation}
  \label{equi}
  u(s,t) = e^{As} z(t)
\end{equation}
where $A\in so(6)$ and $z: \R \lra S^5(1)$, it follows from \eqref{hsphere} that
$u$ is harmonic if and only if $z$ satisfies
\begin{equation}
  \label{neumann1}
  \ddot{z} + A^2z = -( |\dot{z}|^2 + |Az|^2 ) z
\end{equation}
where \ $\dot{ }$ \ denotes differentiation with respect to $t$. As Uhlenbeck noted \cite{U:1982} these
are the equations of motion for the C. Neumann problem of motion of a particle on a 
sphere under the quadratic potential $|Az|^2$.

Define $\R$ actions on $\R^2$ and $S^5(1)$ by
\begin{gather*}
  \gamma \cdot (s,t) = (s+\gamma, t) \\
  \gamma \cdot p = e^{A\gamma}p 
\end{gather*}
where $s,t,\gamma \in \R$ and $p \in S^5(1)$. These induce an action in the usual manner on 
the Banach manifold $C^1(\R^2,S^5)$ by
\begin{equation}
  \label{action}
  \left( \gamma \cdot u \right)(x) = \gamma \cdot u(\gamma^{-1} \cdot x)
\end{equation}
the fixed points of which are precisely maps of the form \eqref{equi}. Since $\R$ acts
by isometries on both $\R^2$ and $S^5$, it follows from the definition of $E$ that it
is an $\R$-invariant function on $C^1(\R^2,S^5)$. Hence we could also appeal to Palais's Principle
of Symmetric Criticality to find the equations satisfied by $z$, as in \cite{U:1982}.

From now on we consider only the case that $A\in u(3)$, so that the one-parameter group $e^{As}$ 
preserves both the metric and the symplectic structure. Then by conjugation we may assume that 
$A=\textrm{diag}\ i(\lambda_1,\lambda_2,\lambda_3).$  In this case equation \eqref{neumann1} becomes
\begin{equation}
  \label{neumann}
  \ddot{z_j} -\lambda_j^2 z_j = -\lambda z_j, \qquad z_j\in \C,\  j=1,2,3
\end{equation}
where 
\begin{equation}
\lambda = |Az|^2 + |\dot{z}|^2.  
\end{equation}

It will be convenient to rewrite the equations slightly.
Writing $z_j=R_j e^{i\theta_j}$ we see that \eqref{neumann} is equivalent to 
\begin{equation}
  \label{radneumann}
  \ddot{R_j} - \frac{J_j^2}{R_j^3} = (\lambda_j^2 - \lambda) R_j\ , \qquad j=1,2,3 
\end{equation}
where $\theta_j$ is determined up to a constant by the relation $J_j=R_j^2\dot{\theta_j}$.

 There are some obvious conserved quantities.
From conservation of energy we have
\begin{align}
  \label{hamiltonian}
&H = |\dot{z}|^2 -|Az|^2 \\
\intertext{and conservation of the quantities}
  \label{angmom}
&J_j = x_j \dot{y_j}-y_j \dot{x_j}\ , \qquad j=1,2,3
\end{align}
expresses the fact that angular momentum in each of the three complex planes $z_1,z_2,z_3$ is conserved.
For details of other less obvious conserved quantities of the Neumann system we refer the reader to 
\cite{U:1982}.

The condition that $u$ be conformal is conveniently expressed in terms of the integrals of motion. 
Namely, $u$ is conformal if and only if
\begin{align}
  \label{conformala}
  &|u_s|^2 - |u_t|^2 = |\dot{z}|^2 - |Az|^2 = H = 0 \\
\intertext{and}
  \label{conformalb}
  &(u_s,u_t) = (\dot{z},Az) = \sum_{i=1}^{3}{\lambda_i J_i}=0.
\end{align}
Summarizing we have
\begin{prop}
\cite{U:1982}
\label{conf_harm}
$u(s,t)=e^{As}z(t):\R^2 \ra S^5(1)$ is a minimal immersion if and only if
$z$ satifies the equations of motion of the Neumann system \eqref{neumann}
and the conserved quantities $H, J_j$ satisfy the constraints \eqref{conformala}
and \eqref{conformalb}.
\end{prop}

\section{$S^1$ Equivariant minimal Legendrian immersions}

\subsection{The Legendrian constraints}
For $u(s,t)=e^{As}z(t)$ to be a minimal Legendrian immersion, besides the conditions
of Proposition \ref{conf_harm}, two further constraints must hold
\begin{align}
  \label{legna}
  &\alpha(u_s) = \omega(u,u_s) = \omega(z,Az) = \sum_{i=1}^3{\lambda_i R_i^2} = 0,\\
  \intertext{and}
  \label{legnb}
  &\alpha(u_t) = \omega(u,u_t) = \omega(z,\dot{z}) = \sum_{i=1}^3{J_i} = 0.
\end{align}
Note that the second equation corresponds merely to further constraints on the values of the integrals
of the Neumann system. The first equation is more mysterious and is in general not
preserved under the flow of the Neumann system. In fact, we have

\begin{lemma}
There are minimal Legendrian immersions of the form given in \eqref{equi} if and only
if $A\in \textrm{su(3)}$.
\end{lemma}

\begin{proof}
From Proposition \ref{minlag_tsleg} any minimal Legendrian immersion
is \tsleg for some $\theta$ and hence $\beta_{\theta}$ restricts to zero on the cone. 
At a point $xu$ on the cone (where $x \in \R^+$) we have
$$ \beta_{\theta}|_{C(u)} = x^2\textrm{Im}\left(e^{i\theta} \det{_{\C}(u, u_s, u_t)} \right) = x^2 \textrm{Im}\left(e^{i\sum{\la_i}s}e^{i\theta}\det{_{\C}(z(t), Az(t), \dot{z}(t))}\right).$$
Since this must hold for all real $s$ and $t$, for $\beta_{\theta}$ to restrict to zero
we must have $\sum{\la_i}=0$ as claimed.

One can also show necessity directly from the equations for a minimal Legendrian equation by showing
that the constraints (\ref{conformala},\ref{conformalb},\ref{legna},\ref{legnb}) are not consistent with the
equations of motion of the Neumann system unless $A \in su(3)$. 

To see this let us compute the second derivative of the mysterious constraint $c:= \omega(z,Az)$ for a solution
of the Neumann system at an instant when all the constraints and their first derivatives are satisfied. One finds
$$  \ddot{c} = \omega(Az,\ddot{z}) + \omega(A\dot{z},\dot{z}) = -\omega(Az,A^2z)+\omega(A\dot{z},\dot{z}).$$
Let $c_1 = \omega(A^2z,Az)$ and $c_2 = \omega(A\dot{z},\dot{z})$.
Then $c_1$ may be expressed in terms of the symmetric polynomials in the $\la_i$ as
$$c_1= \sum{\la_i^3 R_i^2} = (\sum{\la_j})(\sum{\la_i^2 R_i^2})-(\sum{\la_j\la_k})(\sum{\la_i R_i^2})+ \la_1\la_2\la_3 (\sum{R_i^2}).$$
Hence using the constraints we have
$$c_1 =  (\sum{\la_j})\  |Az|^2 + \la_1\la_2\la_3.$$
A calculation shows
$$  c_2 = -\la_1 \la_2 \la_3 \frac{|\dot{z}|^2}{|Az|^2}$$
and so
\begin{equation}
\label{cddot}
  \ddot{c} = (\sum{\la_j})\  |Az|^2 - H \frac{\la_1 \la_2 \la_3}{|Az|^2}.
\end{equation}
Clearly once we have imposed the constraint $H=0$, $\ddot{c} = 0$ if and only if
$A \in su(3)$. Moreover, by differentiating \eqref{cddot} it is easy to verify that 
all higher derivatives of the constraint $c$ also vanish when $A\in su(3)$. Thus to show
existence of minimal Legendrian immersions we need only show there exist initial conditions
for the Neumann system which satisfy all the constraints together with the first derivative of the
mysterious constraint. We will see that this is indeed the case in the proof of Theorem
\ref{min_lag_family} which we now give.
\end{proof}

\begin{proof}[Proof of Theorem \ref{min_lag_family}]
Let $u$ be a minimal Legendrian immersion of the form \eqref{equi}, i.e. 
$u(s,t)=e^{As}z(t)$ where $A\in su(3)$. By conjugation we may assume $A=i\textrm{diag}(\la_1,\la_2,\la_3)$ 
where $\la_1 \geq \la_2 \geq 0 > \la_3$. Let $\alpha={\la_2}/{\la_1}$, then $\alpha \in [0,1]$
and $A = i\textrm{diag}\la_1(1,\alpha,-1-\alpha)$. Moreover, by rescaling $s$ and $t$ we
may assume that $\la_1=1$. Let 
$$\vec{1}=(1,1,1), \; \vec{J}=(J_1,J_2,J_3),\; \vec{\lambda}=(\la_1, \la_2, \la_3), \; \vec{R^2}=(R_1^2, R_2^2, R_3^2).$$
Then the constraints (\ref{conformalb},\ref{legna},\ref{legnb}) together with the 
constraint that $z$ lie on the unit sphere can be written as
\begin{eqnarray}
\label{J_const}
\vec{1}\cdot \vec{J}=0, & \quad & \vec{\la}\cdot \vec{J}=0, \\
\label{R_const}
\vec{1}\cdot \vec{R^2}=0, & \quad &  \vec{\la}\cdot \vec{R^2}=1
\end{eqnarray}
and $A\in su(3)$ is equivalent to $\vec{1}\cdot \vec{\la}=0$. 
Let $\vec{\mu}$ be the cross product of $\vec{1}$ and $\vec{\la}$
$$\vec{\mu}=\vec{1} \times \vec{\la}=(-1-2\alpha, 2 + \alpha, \alpha -1).$$ 
The constraints in \eqref{J_const} are equivalent to 
\begin{equation}
  \label{J}
  \vec{J}=J \vec{\mu}
\end{equation}
for some constant $J$, while the constraints in \eqref{R_const} are equivalent to
\begin{equation}
  \label{gamma}
  \vec{R^2}(t)=\gamma (t) \vec{\mu} + \frac{1}{3}\vec{1}
\end{equation}
for some function $\gamma(t)$. The remaining constraint $|\dot{z}|^2 = |Az|^2$ then becomes
\begin{equation}
  \label{R123}
  \frac{\dot{\gamma}^2}{4} + J^2 = R_1^2R_2^2R_3^2
\end{equation}
or in terms of $\gamma$ 
\begin{equation}
  \label{gamma'}
  \frac{\dot{\gamma}^2}{4} + J^2= \gamma^3 \mu_1\mu_2\mu_3 + \frac{\gamma^2}{3} \sum_{i\neq j}{\mu_i \mu_j} +\frac{1}{27}.
\end{equation}
Since we seek periodic solutions we may assume that $\gamma(0)=\gamma_0>0$, $\dot{\gamma}(0) = 0$. Then at
$t=0$, (\ref{gamma'}) becomes
\begin{equation}
  \label{ga}
  \gamma^3 \mu_1\mu_2\mu_3 + \frac{\gamma^2}{3} \sum_{i\neq j}{\mu_i \mu_j} +\frac{1}{27} = J^2
\end{equation}
and thus specifying $\gamma_0$ determines $J^2$ (and vice versa). 

Let us fix $\alpha \in [0,1]$ and consider the case where $J\neq 0$. 
Given $J \in (0,\frac{1}{3\sqrt{3}}]$, (\ref{ga}) has a unique smallest nonnegative root $\gamma_+(J)$. 
Let $\gamma(0)=\gamma_+(J), \dot{\gamma}(0)=0$. Then up to a translation in time any periodic solution of (\ref{ga}) (except
possibly a solution corresponding to $J=0$ which we shall treat later) arises from such an initial condition. 
Once the initial conditions for $\gamma$ and $\dot{\gamma}$ have been specified,  (\ref{gamma}) fixes 
$R_j(0)$ and $\dot{R}_j(0)$ for $j=1, 2, 3$. Given $J$ and $\alpha$, (\ref{J}) fixes $\vec{J}$. 
If we define $\dot{\theta}_j = J_j/ R_j^2$, then $\dot{\theta_j}(0)$ is determined by (\ref{J}) and (\ref{gamma}) for $j=1,2,3$. By a global rotation in $\textrm{SU(3)}$ (e.g. replacing $z(t)$ by $B z(t)$
where $B=\exp{(i\diag{(\sigma_1, \sigma_2, \sigma_3)})}\in \textrm{SU(3)}$) we may rotate
$z(t)$ so that $\theta_2(0)=\theta_3(0)=0$. We may not assume also that $\theta_1(0)=0$
without allowing $B\in \textrm{U(3)}$ in which case we will change the value of $\theta$
for which $u$ is \tslegg. In the case $J \neq 0$ we shall verify later that choosing 
$-\theta_1(0)=\theta$ gives rise to a \tsleg immersion (in the case $J=0$ choosing 
$-\theta_1(0)= \theta + \pi /2$ gives rise to a \tsleg immersion).

Thus for each $\theta \in [0,2\pi)$, and $(\alpha ,J) \in [0,1] \times (0,1/3\sqrt{3}]$
there is a unique solution of the Neumann equation (given by specifying initial data in the 
manner above) which satisfies the constraints (\ref{conformala},\ref{conformalb},\ref{legna},\ref{legnb}).
Hence by the proof of the 
previous lemma it gives rise to a minimal Legendrian immersion which we denote $u_{\alpha,J}$.

In the case $J=0$ we will explicitly exhibit solutions later in this section and
see that as $\alpha \ra 0$ the period of $\gamma$ becomes infinite, and that the limiting
solution $u_{0,0}$ describes a minimal Legendrian sphere
(which as previously noted is necessarily totally geodesic).

To see which immersions $u_{\alpha,J}$ are geometrically distinct consider in greater detail the geometry of these immersions. 
Since the immersions are all conformal, the metric $g$ induced on $\R^2$ 
can be described by a single positive function $y = |Az|^2 = |\dot{z}|^2$,
where $g=y |dz|^2$. A calculation shows that $\gamma$ and $y$ are related by
\begin{equation}
  \label{az}
  y = -\gamma \mu_1\mu_2\mu_3 + \frac{1}{3}\sum{\la_i^2}.
\end{equation}
It follows from \eqref{az} and \eqref{gamma'} that  $y$ satisfies
\begin{equation}
  \label{ydot}
  \dot{y}^2 + 4y^3 - 2 y^2 \sum{\la_i^2}  = 4C
\end{equation} 
where 
$$-C= \la_1^2\la_2^2\la_3^2 + J^2\mu_1^2\mu_2^2\mu_3^2.$$
The Gauss curvature of the immersion satisfies
\index{Gauss curvature}
\begin{equation}
  \label{gauss}
  K = - \frac{ (\ln y)^{''}}{2y} = 1 + 2C y^{-3}.
\end{equation}

In the case $J=1/3\sqrt{3}$, the corresponding solution of \eqref{gamma'} is $\gamma \equiv 0$ independent of the choice of $\alpha\in [0,1]$ and hence $u_{\alpha,1/3\sqrt{3}}$ 
has $K\equiv 0$. It follows that $u$ must be (a piece of) a generalized Clifford torus. Similarly, if $\alpha =1$ then $\mu_3=0$ and
it follows from \eqref{az} that $y\equiv 2$. Once again $K\equiv 0$ and hence 
$u_{1,J}$ is a (piece of a) generalized Clifford torus.

All other immersions $u_{\alpha,J}$ are geometrically distinct. 
To begin with, note that the remaining immersions are all invariant under a unique $1$-parameter family
of $\textrm{SU(3)}$ -- the subgroup generated by $A=i\diag{(1,\alpha,-1-\alpha)}$.
For $\alpha \in [0,1)$ these are all inequivalent, hence $u_{\alpha,J}$ and $u_{\tilde{\alpha},\tilde{J}}$ are distinct when $\alpha \neq \tilde{\alpha}$.
Now fix $\alpha$ and consider $u_{\alpha, J}$ for $J\in (0,1/3\sqrt{3})$. 
We claim that the minimum and maximum values of the Gauss curvature $K$ are
respectively strictly decreasing and increasing functions of $J$ on $(0,1/3\sqrt{3})$. 
It follows that $u_{\alpha,J}$ and $u_{\alpha,\tilde{J}}$ are geometrically distinct when 
$\alpha \neq \tilde{\alpha}$.

To proof the previous claim, note that \eqref{gauss} shows that for a given immersion
$u_{\alpha,J}$ the minimum (maximum) value of $K$ occurs at the minimum (maximum) value of 
$y$. From (\ref{ydot}) it is clear that for  fixed $\alpha$, $y_{min}$ and $y_{max}$,
the minimum and maximum values attained by y, 
are strictly decreasing and increasing functions of $J$ respectively. 
Since $C$ is a decreasing function of $J$, 
from (\ref{gauss}) we see that the minimum and maximum values of $K$ are, like $y$, strictly decreasing
and increasing functions of $J$ respectively as claimed.
\end{proof}

It is also possible to write down explicit solutions in terms of 
elliptic functions. Let us express $\gamma$ in terms of the Jacobi elliptic functions. 
\index{elliptic functions, Jacobi}
\index{modulus}
Recall that $\gamma$ satisfies the equation
$$   \frac{\dot{\gamma}^2}{4} + J^2= \gamma^3 \mu_1\mu_2\mu_3 + \frac{\gamma^2}{3} \sum_{i\neq j}{\mu_i \mu_j} +\frac{1}{27}$$
and that for $J^2 \in [0,1/27)$ and $\alpha \neq 1$ there are three solutions $\Gamma_1, \Gamma_2, \Gamma_3$ 
to this equation when $\dot{\gamma}=0$. Let us label these solutions so that 
$\Gamma_2 \le 0 \le \Gamma_1 \le \Gamma_3$. Then we can rewrite the previous equation as
\begin{equation}
\label{gam}
  \dot{\gamma}^2 = 4\mu_1 \mu_2 \mu_3 (\gamma - \Gamma_1)(\gamma - \Gamma_2)(\gamma - \Gamma_3).
\end{equation}

\bigskip

\begin{prop}
 $ \gamma(t) = \Gamma_2 - (\Gamma_2 -\Gamma_1) \sn^2{(rt,k)}$ is a solution
of (\ref{gam}) provided
$$ r^2 = \mu_1 \mu_2 \mu_3 (\Gamma_3 - \Gamma_2), \quad \quad  k^2 = \frac{\Gamma_2 - \Gamma_1}{\Gamma_2 - \Gamma_3}$$
where $\sn{}$ is the Jacobi elliptic sn-noidal function.
\end{prop}

\begin{proof}
The proof is a straightforward computation using the basic properties of the Jacobi elliptic functions
(for details see \cite{haskins_thesis}).
\end{proof}
From this proposition and \eqref{gamma} we derive expressions for $R_j^2$
\begin{equation}
  \label{rj}
R_j ^2 = \mu_j (\gamma - \gamma_j) = \mu_j\left((\Gamma_2 - \gamma_j) - (\Gamma_2 - \Gamma_1)\sn^2{(rt,k)}\right)  
\end{equation}
where $\gamma_j=-1/\mu_j$. 

As promised in the proof of Theorem \ref{min_lag_family} we now provide explicit solutions
for the $J=0$ case.

\begin{prop}
For each $\theta\in [0,2\pi)$, 
there exists a family of \tsleg immersions $u_{\alpha,0}: \R^2 \ra S^5(1)$, 
for $\alpha\in [0,1]$, whose Gauss curvature $K$ satisfies \eqref{kmin} and \eqref{kmax} 
(and hence are all distinct). Moreover, $u_{0,0}$ gives rise to a \tsleg sphere
and is the only member of the family $u_{\alpha,J}$ to do so.
\end{prop}

\begin{proof}
In the case $J=0$ we know explicitly  the values of the $\Gamma_i$
$$ \Gamma_i = \gamma_i = -\frac{1}{3\mu_i}, \quad i=1,2,3$$
and hence
\begin{equation}
  \label{r}
 r^2 = (1+2 \alpha), \quad k^2 = \frac{1-\alpha^2}{1+2\alpha}.
\end{equation}
Equation (\ref{rj}) specializes to 
\begin{eqnarray}
 R_1 &=& \mu_1(\gamma_2 - \gamma_1) \cnd{(rt,k)} \\
 R_2 &=& \mu_2(\gamma_1 - \gamma_2) \sn{(rt,k)} \\
 R_3 &=& \mu_3(\gamma_2 - \gamma_3) \dn{(rt,k)}. 
\end{eqnarray}
Define $u_{\alpha,0}$ by the formula
$$u_{\alpha,0}(s,t)= e^{As}(e^{i(\theta+\pi/2)}R_1(t),R_2(t),R_3(t))$$
where as previously we set $A = i\diag{(1,\alpha,-1-\alpha)}$ for $\alpha \in [0,1]$.
Then $u$ is a \tsleg immersion invariant under $e^{As}$.

To find the extreme values taken on by the Gauss curvature, note that
in the case $J=0$ we have
$y_{min} = -\la_2 \la_3 = \alpha (1+ \alpha)$ and $y_{max}=-\la_1\la_3 = 1 + \alpha$.  Thus 
\begin{equation}
  \label{kmin}
  K_{min} = 1 + \frac{2\la_1^2}{\la_2\la_3} = 1 - \frac{2}{\alpha(1+\alpha)}
\end{equation}
and
\begin{equation}
  \label{kmax}
  K_{max}= 1 + \frac{2\la_2^2}{\la_1\la_3} = 1 - \frac{2\alpha^2}{1+\alpha}.
\end{equation}
From \eqref{r}
we see that $k^2 \ra 1$ as $\alpha \ra 0$, and $k^2 \ra 0$ as $\alpha \ra 1$.
In these two limits $\sn$ reduces to $\tanh$ and $\sin$ respectively. 
Thus in the limiting case $\alpha=0, J=0$ we have
$$ \gamma = -\frac{1}{6} + \frac{1}{2}\tanh^2{t}, \quad R_1=R_3=\frac{1}{\sqrt{2}} \sech{t}, \quad R_2 = \tanh{t}.$$

Finally, one can show that in order for any immersion of the form 
$u(s,t)=e^{As}z(t)$ to describe a harmonic sphere, the limit of $z(t)$ as $t\ra \pm \infty$ must be
a fixed point of the action $e^{As}$ \cite{U:1982}. 
Moreover, all the conserved quantities of the Neumann system
must also be zero (since they are zero at a fixed point). For $A\in su(3)$ as above,
$e^{As}$ has nonzero fixed points if and only if $\alpha=0$, in which case any point of the form $(0,z_2,0)\in \C^3$ is fixed.
From equation (\ref{J}), all three angular momenta $J_j$ are zero if and only if $J=0$.
Thus $u_{0,0}$ is the only $u_{\alpha,J}$ which could describe a minimal sphere. In this
case $u$ (in the $0$-\sleg case) has the explicit form
\begin{equation}
\label{sphere}
u(s,t) = (\frac{1}{\sqrt{2}}ie^{is}\sech{t}, \tanh{t},\frac{1}{\sqrt{2}}e^{-is}\sech{t})  
\end{equation}
and we can see directly that the $2$-sphere described is the intersection of the plane 
$$ -i\bar{z_1}=z_3, \quad \textrm{Im}\ z_2 = 0$$
with the $5$-sphere (and hence is totally geodesic). 
\end{proof}

\section{Periodicity conditions}
In order to analyze the periodicity of the immersions $u_{\alpha,J}$ we need the following lemma
whose proof is a short computation (see \cite{haskins_thesis} for full details). 
\begin{lemma}
\label{useful}
For $J\neq 0$ the sum of the angles $\sum{\theta_i}$ and $\dot{\gamma}$ satisfy
\begin{equation}
  \label{theta_sum}
 \dot{\gamma}(t)= 2J \tan{\left(\sum{\theta_i(t)-\theta_i(0)}\right)}.  
\end{equation}
\end{lemma}
Since we chose $\gamma$ so that $\dot{\gamma}(0)=0$, this lemma has the following obvious corollary:
\begin{corollary}
\label{thetasum}
  If $T$ is the period of $\gamma$ then $\sum{\theta_i}(T) = \sum{\theta_i(0)}+n\pi$, for some integer $n$.
\end{corollary}
The previous lemma is also useful in verifying what conditions on $\theta_1(0)$
ensure that the immersions $u_{\alpha,J}$ are \tslegg. For this we need to compute
$\beta_{\theta}$ restricted to the cone on $u$. At a point $(x,s,t)$ on the cone
$$ \beta_{\theta} = x^2\textrm{Im}\left(e^{i\theta} \det{_{\C}(u, u_s, u_t)} \right) = x^2 \textrm{Im}(e^{i\theta}\det{_{\C}(z, Az, \dot{z})}).$$

If $J=0$, so that $\theta_j$ are all constant we have
$$ \det{_{\C}( z,Az,\dot{z})}= i e^{i\sum{\theta_j(0)}} |Az|^2$$
and hence the immersion is \tsleg where $\theta$ depends only on the initial sum
of the angles $\theta_j$. For example, if we choose $\theta_1(0) = \pi/2$ or $\theta_1(0)=3\pi/2$
(and $\theta_2(0)=\theta_3(0)=0$) then the cones are $0$-SLG.

If $J\neq 0$ a short computation using Lemma \ref{useful} shows that
$$ \det{_{\C} (z, Az, \dot{z})} = \frac{i|Az|^2 e^{i\sum{\theta_i(0)}}}{R_1^2R_2^2R_3^2}(J+i\dot{\gamma}/2) (\dot{\gamma}/2+ i J)= -|Az|^2 e^{i\sum{\theta_i(0)}}$$ 
so that now choosing $\theta_1(0)=0$  or $\theta_1(0)=\pi$ gives $0$-\sleg immersions.

\bigskip

Suppose that $(\sigma, \tau)$ is a period of $u(s,t) = e^{As}z(t)$, i.e. 
\begin{equation}
  \label{period}
u(s+\sigma, t+\tau)=u(s,t) \quad \forall s,t.  
\end{equation}
Then the periodicity properties of $u$ are characterized by
\begin{prop} \ \\
\label{period_prop}
 (a) $(\sigma,\tau)$ is a period of $u_{\alpha,J}$ implies $\tau$ is an integer multiple of $T_{\alpha,J}$,
the basic period of $y_{\alpha,J}=|Az|^2$\\
 (b) If $u$ admits two independent periods then it admits a period of the form $(\sigma,0)$\\
 (c) $u$ admits a period of the form $(\sigma,0)$ if and only if $\alpha \in \Q$\\ 
 (d) $u$ admits two independent periods if and only if 

 \begin{equation}
   \label{rationality}
\alpha, \quad  \frac{1}{2\pi}(\alpha \theta_1(T) - \theta_2(T)) \in \Q.
 \end{equation}

\end{prop}

\bigskip

\begin{proof}  
(a) Differentiating (\ref{period}) with respect to $s$ and taking the norm of both sides implies $|Az(t+\tau)|=|Az(t)|$.\\
(b) If the periods are $(\sigma_1, n_1T)$ and $(\sigma_2, n_2T)$ then $(n_1\sigma_2 - n_2\sigma_1,0)$ is also a period.\\
(c) $(\sigma,0)$ is a period implies $e^{i\sigma\la_j} = 1$, for $j=1,2,3$. So $\sigma \la_j \in 2\pi\Z$. In particular, 
$\alpha = \frac{\la_2}{\la_1} \in \Q$. Conversely if $\alpha=\frac{m}{n}$ then $(\frac{2n\pi}{\la_1},0)$ is a period.\\
(d) If $u$ admits two independent periods then $\alpha$ is rational by (b) and (c). By (a) any period is of the form $(\sigma, mT)$. Now since $\dot{\theta_j}$
is $T$-periodic we have $\theta_j(nT)=n\theta_j(T)$. Then periodicity with respect to  $(\sigma, mT)$ is equivalent to 
$e^{i\la_j\sigma + im\theta_j(T)}=1$. Hence $\sigma \la_j + m\theta_j(T) \in 2\pi \Z$, for $j=1,2,3$. Together with 
rationality of $\alpha$ this implies $\alpha \theta_1(T) - \theta_2(T) \in 2\pi \Q$. 

Conversely, by (c) the rationality of
$\alpha$ gives us one period $(\sigma_1,0)$. From above $(\sigma,mT)$ is a period if and only if
$\sigma \la_j + m\theta_j(T) \in 2\pi \Z, \ j=1,2,3$. By assumption 
\mbox{$\frac{1}{2\pi}(\alpha \theta_1(T) - \theta_2(T) = \frac{M}{N}$} for some integers $M$ and $N$.
With $\sigma=-\frac{2N\theta_1(T)}{\la_1}$ and $m=2N$ the period condition becomes
  \begin{displaymath}
    2N\left(-\frac{\la_j}{\la_1} \theta_1(T) + \theta_j(T) \right) \in  2\pi \Z, \quad j=1,2,3.
  \end{displaymath}
For $j=1$ this condition is trivial, while it holds for $j=2$ because
\begin{displaymath}
  2N\left( -\frac{\la_2}{\la_1}\theta_1(T)+\theta_2(T) \right) = -4N\pi \left(\alpha \theta_1(T) - \theta_2(T)\right) = 
-4M\pi
\end{displaymath}
Since $\sum{\la_i}=0$ and by Corollary \ref{thetasum}, $\theta_3(T) = n\pi - \theta_1(T) - \theta_2(T)$, we have
\begin{displaymath}
 2N\left( -\frac{\la_3}{\la_1}\theta_1(T)+\theta_3(T) \right) = 2N\left(\alpha \theta_1(T) - \theta_2(T) + n\pi\right)=
4M\pi + 2\pi Nn.
\end{displaymath}
So the $j=3$ period condition also holds, and hence 
$(-\frac{2N\theta_1(T)}{\la_1}, 2NT)$ is a second period of $u$.
\end{proof}

Two cases of the previous proposition are particularly interesting: when $J=0$ or $\alpha=0$.
For the case $J=0$ we prove the following result which implies Theorem \ref{main_cone_thm} and part (i) of Theorem \ref{tori}.

\begin{prop}
  For $\alpha \in \Q \cap (0,1]$, the immersion $u_{0,\alpha}$ is doubly periodic and
hence gives rise to a minimal Legendrian torus. Further, let $\alpha = \frac{m}{n}$, where $m < n \in \N$
and $(m,n)=1$. If $mn$ is even, then the period lattice of $u_{\alpha,0}$ is rectangular
with basis $\omega_1 = (2n\pi,0)$, $\omega_2 = (0, 4\Ke(k)/r)$. Otherwise the period lattice
is not rectangular and is generated by $\omega_1 = (2n\pi,0)$ and $\omega_3=(n\pi, 2\Ke(k)/r)$.
In either case each such torus $T_{m,n}$ is embedded and its Gauss curvature satisfies
\eqref{kmin} and \eqref{kmax}.
\end{prop}
\no
\textit{Notation in the proposition:} 
$k$ and $r$ are defined as a functions of $\alpha$ by \eqref{r}, and 
$\Ke$ is the complete elliptic integral defined by
$$\Ke(k) = \int_0^{\pi/2} \frac{dx}{\sqrt{1-k^2 \sin^2{x}}}$$
(the period of $\sn{(t,k)}$ is $4\Ke(k)$).

\begin{proof}
Since $J=0$, the $\theta_i$ are constant and the second condition of part (d) of the previous proposition is 
superfluous. Thus the immersion is doubly periodic if and only if $\alpha\in \Q$. Let $\alpha=\frac{m}{n}$.
It is easy to see that $\omega_1=(2n\pi,0)$ and $\omega_2=(0,4\Ke/r)$ belong to the period lattice of $u_{\alpha,0}$. To find the full
period lattice it is sufficient to find all periods in the rectangle $R$ formed by $0$, $\omega_1$, $\omega_2$ and 
$\omega_1 + \omega_2$. Let $(\sigma, \tau)$ be such a period. By part (a) of the previous proposition 
$\tau$ must be a integer multiple of $2\Ke/r$, the basic period of $y$. It is easy to see that
the smallest period of the form $(\sigma,0)$ occurs when $\sigma=2n\pi$ and so we need only deal with
the case $\tau=2\Ke/r$. Using the fact that $\cnd{(t+2\Ke)}=-\cnd{(t)}$, $\sn{(t+2\Ke)}=-\sn{(t)}$, $\dn(t+2\Ke)=\dn(t)$
we find that $(\sigma, 2\Ke/r)$ is a period if and only iff $\sigma$ satisfies
\begin{equation}
  e^{i\sigma} = -1, \quad e^{i\sigma \alpha} = e^{i\sigma m/n} = -1, \quad e^{-i(1+\alpha)}=1.
\end{equation}
Clearly the third equation is implied by the first two. Moreover, the first equation implies
$e^{i m\sigma} = (-1)^m$, whereas the second implies $e^{im \sigma}=(-1)^n$. Hence if either
$m$ or $n$ is even (both cannot be even since we assumed $(m,n)=1$) then these two equations
are inconsistent. Thus there are no further periods and the period lattice is generated by 
$\omega_1$ and $\omega_2$. If both $m$ and $n$ are odd, then one can check that $\sigma=n\pi$
is the unique solution in $[0,2n\pi)$. Hence $w_3= (n\pi, 2\Ke/r)$ is the only new period in the 
rectangle $R$ and in this case the period lattice is generated by $\omega_1$ (or $\omega_2$) and 
$\omega_3$.

Let us show embeddedness in the case where one of $m$, $n$ is even. The other case is similar, but 
a little more involved since the period lattice is not rectangular. We need to show
that if $s, \tilde{s} \in [0,2n\pi)$ and $t, \tilde{t}\in [0,4\Ke/r)$ and $u(s,t)=u(\tilde{s},\tilde{t})$
then $s=\tilde{s}$ and $t=\tilde{t}$. From our explicit formulae for $R_i$ we see that 
$u(s,t)=u(\tilde{s},\tilde{t})$ is equivalent to 
\begin{eqnarray}
  e^{is}\cnd{(t/r)} &= & e^{i\tilde{s}}\cnd{(\tilde{t}/r)}\\
  e^{is}\sn{(t/r)} &= & e^{i\tilde{s}}\sn{(\tilde{t}/r)}\\
  e^{-i(1+\alpha)s}\dn{(t/r)} &= & e^{i(1+\alpha)\tilde{s}}\dn{(\tilde{t}/r)}.
\end{eqnarray}
Certainly this implies $|\cnd{t}/r| = |\cnd{\tilde{t}/r}|$, which implies
there exists some $T\in [0,\Ke]$ such that $rt, r\tilde{t} \in \{T, 2\Ke -T, 2\Ke + T, 4\Ke - T\}$.
If $t$ and $\tilde{t}$ are distinct, there are essentially two different cases, depending on whether $\cnd{t}=-\cnd{\tilde{t}}, \sn{t}=-\sn{\tilde{t}}$
or $\cnd{t}=\pm\cnd{\tilde{t}}, \sn{t}=\mp\sn{\tilde{t}}$. In the first case the three equations above reduce to 
$$e^{i\sigma}= -1, \quad e^{i\alpha\sigma} = -1, \quad e^{-i\sigma(1+\alpha)} = 1$$
where $\sigma = s-\tilde{s}$. 
That is, we have the same equations as occurred in the periodicity part of the proof. 
Since we assumed one of $m$ and $n$ was even, the first two equations are inconsistent
unless $t=\tilde{t}$ in which case $s=\tilde{s}$ is also forced. 
In the second case the equations reduce to 
$$e^{i\sigma}=\pm 1, \quad e^{i\alpha\sigma} = \mp 1, \quad e^{-i\sigma(1+\alpha)} = 1.$$
Clearly, the first two equations are inconsistent with the third one.
\end{proof}

In the case $J\neq 0$, $\alpha=0$, the conditions in part (d) reduce to  $\theta_2(T) \in \Q$. Using the 
explicit expressions given in the previous chapter and properties of elliptic functions one
 can show that viewed as a function of $J$, $\theta_2(T)$ is strictly monotone. Part (ii) of Theorem \ref{tori} follows.

We conclude with the following result which demonstrates the sharpness of the 
pinching results on minimal Legendrian immersions given in parts (ii) and (iii) of
Theorem \ref{yau}.

\begin{theorem}
  \label{almost_flat}
For any $\epsilon >0$ there exists an embedded minimal Legendrian torus $T$ in $S^5$ which is not flat, but
for which $\sup_{x\in T}{|K(x)|} < \epsilon$, where $K$ is the Gauss curvature of $T$. 
\end{theorem}

\begin{proof}
  Consider an immersion $u_{\alpha,J}$ with $J=0$ and $\alpha = 1- \delta$. 
From (\ref{kmin}) and (\ref{kmax}) the minimum and maximum values of the Gauss curvature
are given by
\begin{equation}
  \label{kmin2}
  K_{min}(u_{1-\delta,0}) =- \frac{\delta(3+\delta)}{(1-\delta)(2-\delta)}
\end{equation}
and
\begin{equation}
  \label{kmax2}
  K_{max}(u_{1-\delta,0}) = \frac{\delta (3- 2\delta)}{2-\delta}.
\end{equation}
Certainly for $\delta < \frac{1}{2}$ we have $|K_{min}| < 7 \delta$ and similarly for $K_{max}$.
Since $u_{1-\delta,0}$ gives rise to an embedded minimal Legendrian torus whenever $\alpha \in \Q$,
just choose $\delta \in \Q \cap (0,\epsilon/7)$ and the result is proved.
\end{proof}

\bibliographystyle{amsplain}
\bibliography{paper}
\end{document}